\documentclass{fundam}

\publyear{22}
\papernumber{2102}
\volume{185}
\issue{1}

\usepackage{url} 
\usepackage[ruled,lined]{algorithm2e}
\usepackage{graphicx}
\usepackage{tikz}
\usepackage{arydshln}
\usepackage{amsmath}
\usepackage{graphicx}
\usepackage{amssymb}
\usepackage{picinpar}
\usepackage{algorithmic}
\usepackage{marvosym}
\usepackage{soul}

\usetikzlibrary{automata, positioning, arrows}

\begin{document}

\title{Outer Independent Roman Domination Number of Cartesian Product of Paths and Cycles}

\address{College of Science, Dalian Maritime University, Dalian, China}


\author{Hong Gao\textsuperscript{*}, Daoda Qiu, Shuyan Du, Yiyue Zhao\\
College of Science \\
Dalian Maritime University \\
Dalian, China\\
gaohong{@}dlmu.edu.cn, qiudaoda@dlmu.edu.cn, 2838157641@qq.com, zhaoyy@dlmu.edu.cn
\and Yuansheng Yang\\
School of Computer Science and Technology \\
Dalian University of Technology \\
Dalian, China \\
yangys{@}dlut.edu.cn}

\maketitle

\runninghead{H. Gao et al.}{Outer Independent Roman Domination Number of Cartesian Product of Paths and Cycles}

\begin{abstract}
Given a graph $G$ with vertex set $V$, an outer independent Roman dominating function (OIRDF) is a function $f$ from $V(G)$ to $\{0, 1, 2\}$ for which every vertex with label $0$ under $f$ is adjacent to at least a vertex with label $2$ but not adjacent to another vertex with label $0$. The weight of an OIRDF $f$ is the sum of vertex function values all over the graph, and the minimum of an OIRDF is the outer independent Roman domination number of $G$, denoted as $\gamma_{oiR}(G)$. In this paper, we focus on the outer independent Roman domination number of the Cartesian product of paths and cycles $P_{n}\Box C_{m}$. We determine the exact values of $\gamma_{oiR}(P_n\Box C_m)$ for $n=1,2,3$ and $\gamma_{oiR}(P_n\Box C_3)$ and present an upper bound of $\gamma_{oiR}(P_n\Box C_m)$ for $n\ge 4, m\ge 4$.

\end{abstract}


\section{Introduction}

In this paper, $G=(V, E)$ is a simple, undirected and connected graph with vertex set $V$ and edge set $E$. For a vertex $v\in V$, the open neighbourhood of $v$ is the set $N(v)=\{{u\in V|uv\in E}\} $ and the closed neighbourhood is the set $N[v]=N(v)\cup \{{v}\}$.

Roman domination on graphs is originated from Roman military strategy in which it is required that every location without legions must be adjacent to at least one location with two legions. In 2004, Cockayne et al. \cite{cockayne2004roman} proposed Roman domination on graphs. A Roman dominating function (RDF) on graph $G=(V,E)$ is a function $f: V\rightarrow \{0, 1, 2\}$ such that for each vertex $u\in V$ with $f(u)=0$, there exists at least one vertex $v\in N(u)$ with $f(v)=2$. The weight of a RDF $f$ is $w(f) = \sum_{v\in V} f(v)$. The Roman domination number is the minimum weight of a RDF of $G$, denoted as $\gamma_{R}(G)$. Let $V_i = \{v\in V|f(v)=i\}$ for $i\in\{0,1,2\}$, then a RDF $f$ can be written as $f = (V_0,V_1,V_2)$.

Roman domination is an attractive topic, and many papers have been published on it \cite{luiz2024roman,pavlivc2012roman,poureidi2020algorithmic,liu2021roman}. In order to strengthen defense capability or reduce costs, scholars studied many variations of Roman domination, such as double Roman domination \cite{beeler2016double}, Italian domination \cite{henning2017italian}, triple Roman domination \cite{ahangar2021triple}, perfect Roman domination \cite{yue2020note}, signed Roman domination \cite{abdollahzadeh2014signed}, and so on. 

In Roman domination, if two locations without legions are not adjacent, then the defense capability can be strengthened. In 2017, Ahangar et al. \cite{abdollahzadeh2017outer} introduced the outer independent Roman domination on graphs. Given a graph $G=(V, E)$, an outer independent Roman dominating function (OIRDF) $f=(V_0,V_1,V_2)$ is a RDF and $V_0$ is independent. The weight of an OIRDF is $w(f)=\sum_{v\in V}f(v)=|V_1|+|V_2|$. The outer independent Roman domination number is the minimum weight of an OIRDF of $G$, denoted as $\gamma_{oiR}(G)$. An OIRDF $f$ is a $\gamma_{oiR}$-function if $w(f)=\gamma_{oiR}(G)$.

Scholars are very interested in this topic and have published many related research results.
Mart\'{i}nez et al. \cite{martinez2020outer} obtained some bounds on $\gamma_{oiR}(G)$ in terms of other parameters and provided closed formulas for the outer independent Roman domination number of rooted product graphs.
Poureidi et al. \cite{poureidi2021algorithmic} proposed an algorithm to compute $\gamma _{oiR}(G)$ in ${\mathcal {O}}(|V| )$ time.
Nazari-Moghaddam et al. \cite{nazari2019trees} provided a constructive characterization of trees $T$ whose outer-independent Roman domination number is equal to its Roman domination number.
Dehgardi et al. \cite{dehgardi2021outer} presented an upper bound of $\gamma _{oiR}(T)$ for a tree $T$ of order $n\ge 3$.
Rad et al. \cite{rad2024bounds} presented upper bounds on the outer-independent Roman domination number for unicyclic and bicyclic graphs.

%
%
%

To compute the exact value of the outer independent Roman domination number in Cartesian product graphs is as difficult as to compute its 2-domination number \cite{garzon20222}. In this paper, we study the outer independent Roman domination number of the Cartesian product of paths and cycles $P_{n}\Box C_{m}$. We determine the exact values of $\gamma_{oiR}(P_{n}\Box C_{3})$ and $\gamma_{oiR}(P_{n}\Box C_{m})$ for $n=1, 2, 3$. We also present an upper on $\gamma_{oiR}(P_{n}\Box C_{m})$ for $m, n\ge 4$.

\begin{figure}[h]
\centering
\includegraphics{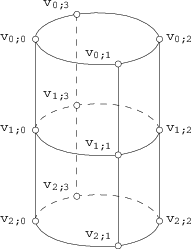}
\caption{Graph $P_{3}\Box C_{4}$.}
\label{fig:p3c4}
\end{figure}

The Cartesian product of paths and cycles is the graph $G=P_{n}\Box C_{m}$ with vertex set $V(G)=\{v_{i,j}|0\leq i\leq n-1, 0\leq j\leq m-1\}$.
$G=P_{n}\Box C_{m}$ is a cylindrical graph and is an important network topology structure graph.
Figure \ref{fig:p3c4} shows graph $P_{3}\Box C_{4}$.
Let $f: V\rightarrow\{0, 1, 2\}$ be an OIRDF on $P_{n}\Box C_{m}$, then we denote $f$ as follows.

\begin{align*}\small{
f=\left(\begin{array}{cccccc}
f(v_{0,0})&f(v_{0,1})&\cdots&f(v_{0,m-1})\\
f(v_{1,0})&f(v_{1,1})&\cdots&f(v_{1,m-1})\\
\vdots    &\vdots    &\ddots&\vdots\\
f(v_{n-1,0})&f(v_{n-1,1})&\cdots&f(v_{n-1,m-1})
\end{array}\right).}
\end{align*}

\section{The outer independent Roman domination number of $P_{1}\Box C_{m}$, $P_{2}\Box C_{m}$ and $P_{3}\Box C_{m}$}

\begin{theorem}\label{thm:p1cm}
For any integer $m\ge 3$, $\gamma_{oiR}(P_{1} \Box C_{m})=3\lfloor\frac{m}{4} \rfloor+t$, where $m\equiv t(\bmod\ 4)$ and $t\in\{0,1,2,3\}$.
\end{theorem}
\begin{proof}
Since $P_{1} \Box C_{m}$ and $C_{m}$ are isomorphic (written $P_{1} \Box C_{m}\cong C_{m}$), then we can achieve the outer independent Roman domination number of $P_{1}\Box C_{m}$ by the outer independent Roman domination number of $C_m$\cite{abdollahzadeh2017outer}.
\end{proof}

\begin{theorem}\label{thm:p2cm}
For any integer $m\ge 3$,
\begin{equation}
\gamma_{oiR}(P_{2} \Box C_{m})=\left\{
\begin{array}{lllllllll}
\lceil\frac{4m}{3}\rceil+1, & m\equiv3,5(\bmod\ 6),\\
\lceil\frac{4m}{3}\rceil,   & otherwise.
\end{array}\notag
\right.
\end{equation}
\end{theorem}
\begin{proof}
Let $G=P_{2}\Box C_{m}$, $V(G)=\{v_{i,j}|0\leq i\leq 1, 0\leq j\leq m-1\}$.

(1) $\gamma_{oiR}(P_{2}\Box C_{m})\le \lceil\frac{4m}{3}\rceil+1$ for $m\equiv3,5(\bmod\ 6)$, $\gamma_{oiR}(P_{2}\Box C_{m})\le \lceil\frac{4m}{3}\rceil$ otherwise.
We first define a function $h$ on $\{v_{i,j}|0\leq i\leq 1, 0\leq j\leq 5\}$.
\begin{eqnarray*}\small{
h(v_{i,j})=\left\{
\begin{array}{llllll}
2,   &(i=0\wedge j=0)\vee (i=1\wedge j=3),\\
1,   &(i=0\wedge j=2,4)\vee (i=1\wedge j=1,5),\\
0,   &otherwise.
\end{array}
\right. \text{\ \ i.e.\ \ }
h=\left(
\begin{array}{llllll}
201010\\
010201
\end{array}
\right).}
\end{eqnarray*}

Then we define an OIRDF $f$ on $P_{2} \Box C_{m}$ as follows,
\begin{equation}\small{
f(v_{i,j})=\left\{
\begin{array}{llllll}
1,                      &m\equiv1(\bmod\ 6), i=0, j=m-2,\\
                        &m\equiv1,3,5(\bmod\ 6), i=1, j=m-1,\\
0,                      &m\equiv1(\bmod\ 6), i=0, j=m-1,\\
h(v_{i,j\bmod6}),   &otherwise.
\end{array}\notag
\right.}
\end{equation}


That is,
\small{
\begin{align*}
m & \equiv 0 \, (\bmod \, 6), & m & \equiv 1 \, (\bmod \, 6), \\
f & = \begin{pmatrix}
201010&\cdots&201010\\
010201&\cdots&010201
\end{pmatrix}, &
f & = \begin{pmatrix}
201010&\cdots&201010& 2010110\\
010201&\cdots&010201& 0102011
\end{pmatrix} \\
m & \equiv 2 \, (\bmod \, 6), & m & \equiv 3 \, (\bmod \, 6), \\
f & = \begin{pmatrix}
201010&\cdots&201010&20\\
010201&\cdots&010201&01
\end{pmatrix}, &
f & = \begin{pmatrix}
201010&\cdots&201010&201\\
010201&\cdots&010201&011
\end{pmatrix} \\
m & \equiv 4 \, (\bmod \, 6), & m & \equiv 5 \, (\bmod \, 6), \\
f & = \begin{pmatrix}
201010&\cdots&201010&2010\\
010201&\cdots&010201&0102
\end{pmatrix}, &
f & = \begin{pmatrix}
201010&\cdots&201010&20101\\
010201&\cdots&010201&01021
\end{pmatrix}
\end{align*}}
One can check that $f$ is an OIRDF and the weight of $f$ is

\begin{equation}\small{
w(f)=\left\{
\begin{array}{lllllllll}
\frac{m}{6}\times8=\frac{4m}{3},     & m\equiv0(\bmod\ 6),\\
\frac{m-7}{6}\times8+10=\frac{4m+2}{3},     & m\equiv1(\bmod\ 6),\\
\frac{m-2}{6}\times8+3=\frac{4m+1}{3},     & m\equiv2(\bmod\ 6),\\
\frac{m-3}{6}\times8+5=\frac{4m+3}{3},     & m\equiv3(\bmod\ 6),\\
\frac{m-4}{6}\times8+6=\frac{4m+2}{3},     & m\equiv4(\bmod\ 6),\\
\frac{m-5}{6}\times8+8=\frac{4m+4}{3},     & m\equiv5(\bmod\ 6).\\
\end{array}\notag
\right.}
\end{equation}

Hence, we have
\begin{equation}
\gamma_{oiR}(P_{2} \Box C_{m})\le \left\{
\begin{array}{lllllllll}
\lceil\frac{4m}{3}\rceil+1, & m\equiv3,5(\bmod\ 6),\\
\lceil\frac{4m}{3}\rceil,  & otherwise.
\end{array}\notag
\right.
\end{equation}

(2) $\gamma_{oiR}(P_{2}\Box C_{m})\ge \lceil\frac{4m}{3}\rceil+1$ for $m\equiv3,5(\bmod\ 6)$, $\gamma_{oiR}(P_{2}\Box C_{m})\ge \lceil\frac{4m}{3}\rceil$ otherwise.

Let $f$ be a $\gamma_{oiR}$-function of $P_{2}\Box C_m$. Denote $V^{j}=\{v_{0,j}, v_{1,j}\}(0\leq j \leq m-1)$, $f^{j}=f(V^{j})=f(v_{0,j})+f(v_{1,j})$.
By the definition of OIRDF, $f^{j}\ge 1$ and if $f^{j}=1$, then $f^{j-1}\ge 2$ or $f^{j+1}\ge 2$.
Then we use the following algorithm to group $V^{j}$ ($0\le j\le m-1$) into three categories.

\begin{algorithm}[h!] \label{bagging:p2cm}\caption{Dividing groups $V^j$ ($0\le j\le m-1,j\in N$) into three categories}
\KwIn{$V^j$}
\KwOut{$t_0,t_1,t_2,B^{0}_{t_0},B^{1}_{t_1},B^{2}_{t_2}$}
\SetKw{Init}{Initialize}
\Init{$t_0=t_1=t_2=0$, $B^{0}_{t_0}=B^{1}_{t_1}=B^{2}_{t_2}=\emptyset$ and $d^{j}=0$ for $j=0$ up to $m-1$}\;
\For{$j$ from $0$ to $m-1$ with $f^{j}\geq 3\wedge d^{j}=0$}
{
$t_2=t_2+1$, $d^{j}=1$, $B^{2}_{t_2}=\{V^{j}\}$;

  \If{$f^{j-1}=1\wedge d^{j-1}=0$}
  {
    $d^{j-1}=1$, $B^{2}_{t_2}=B^{2}_{t_2}\cup \{V^{j-1}\}$\;
  }
  \If{$f^{j+1}=1\wedge d^{j+1}=0$}
  {
    $d^{j+1}=1$, $B^{2}_{t_2}=B^{2}_{t_2}\cup \{V^{j+1}\}$\;
  }
  $Note: \sum_{V^{j}\in B^{2}_{t_2}}f^{j}> \frac{4}{3}|B^{2}_{t_2}|+\frac{2}{3}$.
}

\For{$j$ from $0$ to $m-1$ with $f^{j}=2\wedge d^{j}=0$}
{
  \If{$f^{j-1}=f^{j+1}=1\wedge d^{j-1}=d^{j+1}=0$}
  {
    $t_{0}=t_{0}+1$, $d^{j-1}=d^{j}=d^{j+1}=1$, $B^{0}_{t_0}=\{V^{j-1},V^{j},V^{j+1}\}$\;
    $Note: \sum_{V^{j}\in B^{0}_{t_0}}f^{j}=4=\frac{4}{3}|B^{0}_{t_0}|$.
  }
  \ElseIf{$f^{j-1}=1\wedge d^{j-1}=0$}
  {
    $t_{1}=t_{1}+1$, $d^{j-1}=d^{j}=1$, $B^{1}_{t_1}=\{V^{j-1},V^{j}\}$\;
    $Note: \sum_{V^{j}\in B^{1}_{t_1}}f^{j}=3=\frac{4}{3}|B^{1}_{t_1}|+\frac{1}{3}$.
  }
  \ElseIf{$f^{j+1}=1\wedge d^{j+1}=0$}
  {
    $t_{1}=t_{1}+1$, $d^{j}=d^{j+1}=1$, $B^{1}_{t_1}=\{V^{j},V^{j+1}\}$\;
    $Note: \sum_{V^{j}\in B^{1}_{t_1}}f^{j}=3=\frac{4}{3}|B^{1}_{t_1}|+\frac{1}{3}$.
  }
  \Else
  {
    $t_{2}=t_{2}+1$, $d^{j}=1$, $B^{2}_{t_2}=\{V^{j}\}$\;
    $Note: \sum_{V^{j}\in B^{2}_{t_2}}f^{j}=2=\frac{4}{3}|B^{2}_{t_2}|+\frac{2}{3}$.
  }
}
\end{algorithm}

By Algorithm \ref{bagging:p2cm},
\begin{equation*}
\begin{split}
w(f)&=\sum\limits_{j=0}^{m-1}f^j\\
    &=\sum\limits_{t=1}^{t_0}\sum_{V^j\in B^{0}_{t}}f^j+\sum\limits_{t=1}^{t_1}\sum_{V^j\in B^{1}_{t}}f^j+\sum\limits_{t=1}^{t_2}\sum_{V^j\in B^{2}_{t}}f^j\\
    &\geq\sum\limits_{t=1}^{t_0}\frac{4}{3}|B^{0}_{t}|+\sum\limits_{t=1}^{t_1}(\frac{4}{3}|B^{1}_{t}|+\frac{1}{3})+\sum\limits_{t=1}^{t_2}(\frac{4}{3}|B^{2}_{t}|+\frac{2}{3})\\
    &=\frac{4m}{3}+\frac{t_1}{3}+\frac{2t_2}{3}.
\end{split}
\end{equation*}

Hence, $w(f)\ge \lceil\frac{4m}{3}\rceil$. Then $\gamma_{oiR}(P_{2}\Box C_{m})\ge \lceil\frac{4m}{3}\rceil$ for $m\not\equiv3,5(\bmod\ 6)$. Next, we will prove $\gamma_{oiR}(P_{2}\Box C_{m})\ge \lceil\frac{4m}{3}\rceil+1$ for $m\equiv3,5(\bmod\ 6)$ by contradiction.

Suppose $w(f)\le \lceil\frac{4m}{3}\rceil$ for $m\equiv3,5(\bmod\ 6)$.
For $m\equiv3(\bmod\ 6)$, $w(f)\le \lceil\frac{4m}{3}\rceil=\frac{4m}{3}$, then $t_1=t_2=0$.
For $m\equiv5(\bmod\ 6)$, $w(f)\le \lceil\frac{4m}{3}\rceil=\frac{4m}{3}+\frac{1}{3}$, then $t_1\le 1$, $t_2=0$.
By Algorithm \ref{bagging:p2cm}, we know the details on $B^{k}_{t}$ ($0\le k\le 1$, $0\le t\le t_k$).

$B^{0}_{t}=\{V^{j-1},V^{j},V^{j+1}\}$ ($0\le t\le t_0$) where $f(V^{j})=2$ and $f(V^{j-1})=f(V^{j+1})=1$, then we denote
$f(B^{0}_{t})=(1,2,1)$, and the corresponding segments of $f$ are as follows.
$$\small{\begin{array}{ccccccccccccccccc}
S_{0,1}=\left(
\begin{array}{ccccccc}
010\\
111
\end{array}\right),
S_{0,2}=\left(
\begin{array}{ccccccc}
111\\
010
\end{array}\right),
S_{0,3}=\left(
\begin{array}{ccccccc}
110\\
011
\end{array}\right),\\
S_{0,4}=\left(
\begin{array}{ccccccc}
011\\
110
\end{array}\right),
S_{0,5}=\left(
\begin{array}{ccccccc}
020\\
101
\end{array}\right),
S_{0,6}=\left(
\begin{array}{ccccccc}
101\\
020
\end{array}\right).
\end{array}}$$

$B^{1}_{t}=\{V^{j-1},V^{j}\}$ or $\{V^{j},V^{j+1}\}$ ($0\le t\le t_1$) where $f(V^{j})=2$ and $f(V^{j-1})=f(V^{j+1})=1$, then
$f(B^{1}_{t})=(1,2)$ or $(2,1)$, and the corresponding segments of $f$ are as follows.
$$\small{\begin{array}{ccccccccccccccccc}
S_{1,1}=\left(
\begin{array}{ccccccc}
01\\
11
\end{array}\right),
S_{1,2}=\left(
\begin{array}{ccccccc}
11\\
01
\end{array}\right),
S_{1,3}=\left(
\begin{array}{ccccccc}
02\\
10
\end{array}\right),
S_{1,4}=\left(
\begin{array}{ccccccc}
10\\
02
\end{array}\right),\\
S_{1,5}=\left(
\begin{array}{ccccccc}
10\\
11
\end{array}\right),
S_{1,6}=\left(
\begin{array}{ccccccc}
11\\
10
\end{array}\right),
S_{1,7}=\left(
\begin{array}{ccccccc}
20\\
01
\end{array}\right),
S_{1,8}=\left(
\begin{array}{ccccccc}
01\\
20
\end{array}\right).
\end{array}}$$

For $m\equiv3(\bmod\ 6)$, $t_1=t_2=0$. Then $f$ \textbf{only} has $S_{0,i}$ ($1\le i\le 6$).
In fact, $f$ \textbf{only} has $S_{0,5}$ and $S_{0,6}$.
If $f$ contains $S_{0,1}$, then
$$\small{
f(V)=\left (\begin{array}{ccccccc}
\cdots&f(v)&010&f(u)&\cdots\\
\cdots&    &111&    &\cdots
\end{array}\right )}$$
where $f(v)=f(u)=2$. It has $S_{0,1}$ cannot be placed before or after $S_{0,i}$ ($1\le i\le 6$). Thus, $f$ cannot contains $S_{0,1}$. Similarly, $f$ can not contains $S_{0,i}$ ($2\le i\le 4$).
Since $f$ is an OIRDF, it cannot consist of only $S_{0,5}$ or only $S_{0,6}$. Then it must be
$$\small{
f(V)=\left (\begin{array}{cc:c:cc:cc}
020&101&\cdots&020&101&020&101  \\
101&020&\cdots&101&020&101&020
              \end{array}
           \right ).}
$$
It has $m\equiv0(\bmod\ 6)$ which contradicts to $m\equiv3(\bmod\ 6)$.

For $m\equiv5(\bmod\ 6)$, $t_1\le 1$, $t_2=0$.
If $t_1=0$, then $m\equiv0(\bmod\ 6)$, thus $t_1=1$ for $m\equiv5(\bmod\ 6)$. That is, $f$ contains \textbf{only one} $S_{1,i}$ $(1\le i\le 8)$. Since $f$ is an OIRDF, it cannot contain $S_{1,1}$, $S_{1,2}$, $S_{1,5}$ and $S_{1,6}$. Then it must be
$$\begin{array}{lllllllllll}
f(V)=\small{\left (\begin{array}{cc:c:ccc:c:cc}
020&101&\cdots&020&10&101&\cdots&020&101  \\
101&020&\cdots&101&02&020&\cdots&101&020
              \end{array}
           \right )}\mbox{ or }\\
f(V)=\small{\left (\begin{array}{cc:c:ccc:c:cc}
020&101&\cdots&020&20&101&\cdots&020&101  \\
101&020&\cdots&101&01&020&\cdots&101&020
              \end{array}
           \right )}\mbox{ or }\\
f(V)=\small{\left (\begin{array}{cc:c:ccc:c:cc}
101&020&\cdots&101&02&020&\cdots&101&020  \\
020&101&\cdots&020&10&101&\cdots&020&101
              \end{array}
           \right )}\mbox{ or }\\
f(V)=\small{\left (\begin{array}{cc:c:ccc:c:cc}
101&020&\cdots&101&01&020&\cdots&101&020  \\
020&101&\cdots&020&20&101&\cdots&020&101
              \end{array}
           \right )}
\end{array}$$
It has $m\equiv2(\bmod\ 6)$ which contradicts to $m\equiv5(\bmod\ 6)$.

Therefore, $\gamma_{oiR}(P_{2}\Box C_{m})=w(f)\ge \lceil\frac{4m}{3}\rceil+1$ for $m\equiv3,5(\bmod\ 6)$.
\end{proof}

\begin{theorem}\label{thm:p3cm}
For any integer $m\ge 3$,
\begin{equation}
\gamma_{oiR}(P_{3} \Box C_{m})=\left\{
\begin{array}{lllllllll}
2m,       & m\equiv0(\bmod\ 2),\\
2m+1,     & m\equiv1(\bmod\ 2).
\end{array}\notag
\right.
\end{equation}
\end{theorem}
\begin{proof}
Let $G=P_{3} \Box C_{m}$, $V(G)=\{v_{i,j}|0\leq i\leq 2, 0\leq j\leq m-1\}$.

(1) $2m$ and $2m+1$ are the upper bounds of $\gamma_{oiR}(P_{3} \Box C_{m})$ for $m\equiv0(\bmod\ 2)$ and $m\equiv1(\bmod\ 2)$ respectively. We first define a function $h$ on $\{v_{i,j}|0\leq i\leq 2, 0\leq j\leq 1\}$.

\begin{eqnarray*}\small{
h(v_{i,j})=\left\{
\begin{array}{llllll}
2,   &i=1, j=0,\\
1,   &i=0,2, j=1,\\
0,   &otherwise.
\end{array}
\right. \text{\ \ i.e.\ \ }
h=\left(
\begin{array}{llllll}
01\\
20\\
01
\end{array}
\right).}
\end{eqnarray*}

Then we define an OIRDF $f$ on $P_{3} \Box C_{m}$ as follows,
\begin{equation}\small{
f(v_{i,j})=\left\{
\begin{array}{llllll}
1,                    &m\equiv1(\bmod\ 2), j=m-1,\\
h(v_{i,j\bmod\ 2}),   &otherwise.
\end{array}\notag
\right.}
\end{equation}

That is,
\begin{align*}\small{
\begin{array}{llllllllll}
m\equiv0(\bmod\ 2),&&&&&&m\equiv1(\bmod\ 2),\\
f=\left(\begin{array}{ccccccccccc}
01&\cdots&01\\
20&\cdots&20\\
01&\cdots&01
\end{array}
\right),&&&&&&
f=\left(\begin{array}{cccccc}
01&\cdots&01&1\\
20&\cdots&20&1\\
01&\cdots&01&1
\end{array}\right).
\end{array}}
\end{align*}

One can check $f$ is an OIRDF and the weight is $w(f)=4\times\frac{m}{2}=2m$ for $m\equiv0(\bmod\ 2)$ and $w(f)=4\times\frac{m-1}{2}+3=2m+1$ for $m\equiv1(\bmod\ 2)$.

Hence, $\gamma_{oiR}(P_{3} \Box C_{m})\le 2m$ for $m\equiv0(\bmod\ 2)$ and $\gamma_{oiR}(P_{3} \Box C_{m})\le 2m+1$ for $m\equiv1(\bmod\ 2)$.

(2) $2m$ and $2m+1$ are the lower bounds of $\gamma_{oiR}(P_{3} \Box C_{m})$ for $m\equiv0(\bmod\ 2)$ and $m\equiv1(\bmod\ 2)$ respectively.

We can find a $\gamma_{oiR}$-function $f$ satisfying properties (a) and (b). Denote $V^{j}=\{v_{i,j}\mid 0\leq i\leq 2\}(0\leq j \leq m-1)$, $f^{j}=f(V^{j})=\sum_{v_{i,j}\in V^{j}}f(v_{i,j})$.

(a) For $0\leq j \leq m-1$, $f^{j}\ge1$; if $f^{j}=1$, then $f^{j+1}\ge 2$ and $f^{j-1}+f^{j+1}\geq 6$.

(b) If $f^{j-1}=f^{j+1}=1$, then $2\le f^{j}\le3$.

\noindent In which superscripts are taken modulo $m$.

For (a), since $V_0$ is independent, then $f^{j}\neq 0$, i.e. $f^{j}\ge1$.
Furthermore, if $f^{j}=1$, then $f(v_{0,j})=f(v_{2,j})=0$, $f(v_{1,j})=1$, and $f(v_{0,j+1})\ge 1$ and $f(v_{2,j+1})\ge 1$, it follows $f^{j+1}\ge 2\neq 1$. Additionally, since for each $u\in V_0$, there exists at least one vertex $v\in N(u)$ with $f(v)=2$, then $f(v_{0,j-1})+f(v_{0,j+1})\geq3$, $f(v_{2,j-1})+f(v_{2,j+1})\geq3$, it follows $f^{j-1}+f^{j+1}\geq 6$.

For (b), if $f^{j-1}=f^{j+1}=1$, by (a), $f^{j}\neq 1$. Suppose $f^{j}\ge 4$. Since $V_{0}$ is independent and $f^{j-1}=f^{j+1}=1$, then
$f(v_{0,j-1})=f(v_{2,j-1})=f(v_{0,j+1})=f(v_{2,j+1})=0$, $f(v_{1,j-1})=f(v_{1,j+1})=1$.
Define an OIRDF $f'$ as: $f'(v_{0,j-1})=f'(v_{2,j-1})=f'(v_{1,j})=f'(v_{0,j+1})=f'(v_{2,j+1})=0$, $f'(v_{0,j})=f'(v_{2,j})=1$, $f'(v_{1,j-1})=f'(v_{1,j+1})=2$.
Then $w(f')\le w(f)$. If $w(f')< w(f)$, then there is contradiction to $f$ being a $\gamma_{oiR}$-function. If $w(f')=w(f)$, then we find a $\gamma_{oiR}$-function $f'$ with $f'_{j-1}=f'_{j+1}=1$ and $2\le f'_{j}\le3$.

Then we use Algorithm \ref{bagging:p3cm} to group $V^{j}$ ($0\le j\le m-1$) into two categories.

\begin{algorithm}[h!] \label{bagging:p3cm}\caption{Dividing groups $V^j$ ($0\le j\le m-1,j\in N$) into two categories}
\KwIn{$V^j$}
\KwOut{$t_0,t_1,B^{0}_{t_0},B^{1}_{t_1}$}
\SetKw{Init}{Initialize}
\Init{$t_0=t_1=0$, $B^{0}_{t_0}=B^{1}_{t_1}=\emptyset$ and $d^{j}=0$ for $j=0$ up to $m-1$}\;
\For{$j$ from 0 to $m-1$ with $f^{j}\geq 4\wedge d^{j}=0$}
{
$t_1=t_1+1$, $d^{j}=1$, $B^{1}_{t_1}=\{V^{j}\}$;

  \If{$f^{j-1}=1\wedge d^{j-1}=0$}
  {
    $d^{j-1}=1$, $B^{1}_{t_1}=B^{1}_{t_1}\cup \{V^{j-1}\}$\;
  }
  \If{$f^{j+1}=1\wedge d^{j+1}=0$}
  {
    $d^{j+1}=1$, $B^{1}_{t_1}=B^{1}_{t_1}\cup \{V^{j+1}\}$\;
  }
  $Note: \sum_{V^{j}\in B^{1}_{t_1}}f^{j}\ge 2|B^{1}_{t_1}|+1$.
}

\For{$j$ from 0 to $m-1$ with $f^{j}=3\wedge d^{j}=0$}
{
  \If{$f^{j+1}=1\wedge d^{j+1}=0$}
  {
    $t_{0}=t_{0}+1$, $d^{j}=d^{j+1}=1$, $B^{0}_{t_0}=\{V^{j},V^{j+1}\}$\;
    $Note: \sum_{V^{j}\in B^{0}_{t_0}}f^{j}=4=2|B^{0}_{t_0}|$.
  }
  \Else
  {
    $t_{1}=t_{1}+1$, $d^{j}=1$, $B^{1}_{t_1}=\{V^{j}\}$\;
    $Note: \sum_{V^{j}\in B^{1}_{t_1}}f^{j}=3=2|B^{1}_{t_1}|+1$.
  }
}

\For{$j$ from 0 to $m-1$ with $f^{j}=2\wedge d^{j}=0$}
{
  $t_0=t_0+1$, $d^{j}=1$, $B^{0}_{t_0}=\{V^{j}\}$;

  $Note: \sum_{V^{j}\in B^{0}_{t_0}}f^{j}=2=2|B^{0}_{t_0}|$.

}
\end{algorithm}

By Algorithm \ref{bagging:p3cm},
\begin{equation*}
\begin{split}
w(f)&=\sum\limits_{j=0}^{m-1}f^j=\sum\limits_{t=1}^{t_0}\sum_{V^j\in B^{0}_{t}}f^j+\sum\limits_{t=1}^{t_1}\sum_{V^j\in B^{1}_{t}}f^j
    \geq\sum\limits_{t=1}^{t_0}(2|B^{0}_{t}|)+\sum\limits_{t=1}^{t_1}(2|B^{1}_{t}|+1)
    =2m+t_1.
\end{split}
\end{equation*}

Hence, $w(f)\ge 2m$. Then $\gamma_{oiR}(P_{3}\Box C_{m})\ge 2m$ for $m\equiv0(\bmod\ 2)$. Next, we will prove $\gamma_{oiR}(P_{3}\Box C_{m})\ge 2m+1$ for $m\equiv1(\bmod\ 2)$ by contradiction.

Suppose $w(f)\le 2m$ for $m\equiv1(\bmod\ 2)$, then $t_1=0$. By Algorithm \ref{bagging:p3cm}, we know the details on $B^{0}_{t}$ ($0\le t\le t_0$).

$B^{0}_{t}=\{V^{j},V^{j+1}\}$ where $f(V^{j})=3$ and $f(V^{j+1})=1$, then $f(B^{0}_{t})=(3,1)$.
Or $B^{0}_{t}=\{V^{j}\}$ where $f(V^{j})=2$, then $f(B^{0}_{t})=(2)$.
The corresponding segments of $f$ are as follows.
$$\small{\begin{array}{lllllllllccccccccc}
&S_{0,1}=\left(
\begin{array}{ccccccc}
10\\
11\\
10
\end{array}\right),
&S_{0,2}=\left(
\begin{array}{ccccccc}
10\\
01\\
20
\end{array}\right),
&S_{0,3}=\left(
\begin{array}{ccccccc}
20\\
01\\
10
\end{array}\right),\\
&S_{0,4}=\left(
\begin{array}{ccccccc}
0\\
1\\
1
\end{array}\right),
&S_{0,5}=\left(
\begin{array}{ccccccc}
1\\
1\\
0
\end{array}\right),
&S_{0,6}=\left(
\begin{array}{ccccccc}
1\\
0\\
1
\end{array}\right),
&S_{0,7}=\left(
\begin{array}{ccccccc}
0\\
2\\
0
\end{array}\right).
\end{array}}$$

Since $t_1=0$, then $f$ \textbf{only} contains $S_{0,i}$ $(1\le i\le 7)$.
However, $f$ does not contain $S_{0,1}$. If not, then
$$\small{
f(V)=\left (\begin{array}{ccccccc}
\cdots&10&f(u)&\cdots\\
\cdots&11&    &\cdots\\
\cdots&10&f(v)&\cdots
\end{array}\right )}$$
where $f(u)=f(v)=2$. It has $S_{0,1}$ cannot be placed before $S_{0,i}$ ($1\le i\le 7$). Therefore, $f$ does not contain $S_{0,1}$.

$f$ does not contain $S_{0,2}$ and $S_{0,3}$, otherwise $f$ must be
$$\small{
f(V)=\left (\begin{array}{ccccccc}
10&20&\cdots&10&20\\
01&01&\cdots&01&01\\
20&10&\cdots&20&10
\end{array}\right )}.$$
It has $m\equiv0(\bmod\ 2)$ which contradicts to $m\equiv1(\bmod\ 2)$.

$f$ does not contain $S_{0,4}$ and $S_{0,5}$. If $f$ contains $S_{0,4}$, then
$$\small{
f(V)=\left (\begin{array}{ccccccc}
\cdots&f(u)&0&f(v)&\cdots\\
\cdots&    &1&    &\cdots\\
\cdots&    &1&    &\cdots
\end{array}\right )}$$
where $f(u)=2$ or $f(v)=2$. It follows $S_{0,4}$ cannot be placed before or after $S_{0,i}$ ($4\le i\le 7$). Since $f$ does not contain $S_{0,i}$ ($i=1,2,3$), then $f$ does not contain $S_{0,4}$. Similarly, $f$ does not contain $S_{0,5}$.

Then $f$ consists of $S_{0,6}$ and $S_{0,7}$. Since $f$ is an OIRDF, then $f$ must be
$$\small{
f(V)=\left (\begin{array}{ccccccc}
10&10&\cdots&10&10\\
02&02&\cdots&02&02\\
10&10&\cdots&10&10
\end{array}\right )}.$$
It has $m\equiv0(\bmod\ 2)$ which contradicts to $m\equiv1(\bmod\ 2)$.

Therefore, $\gamma_{oiR}(P_{3}\Box C_{m})=w(f)\ge 2m+1$ for $m\equiv1(\bmod\ 2)$.
\end{proof}

\section{The outer independent Roman domination number of $P_{n}\Box C_{3}$}

\begin{theorem}\label{thm:pnc3}
For any integer $n\ge 3$, $\gamma_{oiR}(P_{n} \Box C_{3})=\lceil\frac{7n}{3}\rceil$.
\end{theorem}
\begin{proof}
Let $G=P_{n} \Box C_{3}$, $V(G)=\{v_{i,j}|0\leq i\leq n-1, 0\leq j\leq 2\}$.

(1) $\lceil\frac{7n}{3}\rceil$ is the upper bound of $\gamma_{oiR}(P_{n} \Box C_{3})$. We first define a function $h$ on $P_{6} \Box C_{3}$.
\begin{equation}\small{
h(v_{i,j})=\left\{
\begin{array}{llllll}
2,   &(i=1, j=0)\vee (i=4, j=1),\\
1,   &(i=3,5, j=0)\vee (i=0,2, j=1)\vee j=2,\\
0,   &otherwise.
\end{array}\notag
\right. \text{\ \ i.e.\ \ }
h=\left(
\begin{array}{llllll}
011\\
201\\
011\\
101\\
021\\
101
\end{array}
\right).}
\end{equation}

Then we define an OIRDF $f$ on $P_{3} \Box C_{m}$ as follows.
\begin{equation}\small{
f(v_{i,j})=\left\{
\begin{array}{llllll}
2,   &n\equiv1,4(\bmod\ 6), i=n-1, j=2,\\
h(v_{i\bmod6,j}),   &otherwise.
\end{array}\notag
\right.}
\end{equation}

That is,
\begin{eqnarray*}\small{
\left(
\begin{array}{cc}
011\\
201\\
011\\
101\\
021\\
101\\
\vdots\\
011\\
201\\
011\\
101\\
021\\
101\\
\end{array}
\right),
\left(
\begin{array}{cc}
011\\
201\\
011\\
101\\
021\\
101\\
\vdots\\
011\\
201\\
011\\
101\\
021\\
101\\
\\
012
\end{array}
\right),
\left(
\begin{array}{cc}
011\\
201\\
011\\
101\\
021\\
101\\
\vdots\\
011\\
201\\
011\\
101\\
021\\
101\\
\\
011\\
201\\
\end{array}
\right),
\left(
\begin{array}{cc}
011\\
201\\
011\\
101\\
021\\
101\\
\vdots\\
011\\
201\\
011\\
101\\
021\\
101\\
\\
011\\
201\\
011\\
\end{array}
\right),
\left(
\begin{array}{cc}
011\\
201\\
011\\
101\\
021\\
101\\
\vdots\\
011\\
201\\
011\\
101\\
021\\
101\\
\\
011\\
201\\
011\\
102\\
\end{array}
\right),
\left(
\begin{array}{cc}
011\\
201\\
011\\
101\\
021\\
101\\
\vdots\\
011\\
201\\
011\\
101\\
021\\
101\\
\ \\
011\\
201\\
011\\
101\\
021\\
\end{array}
\right).
}
\end{eqnarray*}

One can check that $f$ is an OIRDF and the weight is
\begin{equation}\small{
w(f)=\left\{
\begin{array}{lllllllll}
\frac{n}{6}\times14=\frac{7n}{3},     & n\equiv0(\bmod6),\\
\frac{n-1}{6}\times14+3=\frac{7n+2}{3},     & n\equiv1(\bmod6),\\
\frac{n-2}{6}\times14+5=\frac{7n+1}{3},     & n\equiv2(\bmod6),\\
\frac{n-3}{6}\times14+7=\frac{7n}{3},       & n\equiv3(\bmod6),\\
\frac{n-4}{6}\times14+10=\frac{7n+2}{3},       & n\equiv4(\bmod6),\\
\frac{n-5}{6}\times14+12=\frac{7n+1}{3},       & n\equiv5(\bmod6).
\end{array}\notag
\right.}
\end{equation}

Hence, $\gamma_{oiR}(P_{n} \Box C_{3})\le \lceil\frac{7n}{3}\rceil$.

(2) $\lceil\frac{7n}{3}\rceil$ is the lower bound of $\gamma_{oiR}(P_{n} \Box C_{3})$. Let $f$ be a $\gamma_{oiR}$-function of $P_{n} \Box C_{3}$. Denote $V^{i}=\{v_{i,j}\mid 0\leq j\leq 2\}(0\leq i \leq n-1)$, $f^{i}=f(V^{i})=\sum_{v_{i,j}\in V^{i}}f(v_{i,j})$. Then $f$ has properties (a) and (b).

(a) $f^{i-1}+f^i+f^{i+1}\ge 7$ for every integer $i$ $(1\leq i\leq n-2)$ where superscripts are taken modulo $n$.

Since $V_0$ is independent, we have $f^i\ge 2$ for every $i$ $(0\le i\le n-1)$. For some $i$ $(1\le i\le n-2)$,
if $f^{i}\ge 3$, then $f^{i-1}+f^{i}+f^{i+1}\ge 7$;
if $f^{i}=2$, without loss of generality, let $f(v_{i,0})=0$, $f(v_{i,1})=f(v_{i,2})=1$,
then $f(v_{i-1,0})+f(v_{i+1,0})\ge 3$, $f(v_{i-1,1})+f(v_{i-1,2})\ge 1$, $f(v_{i+1,1})+f(v_{i+1,2})\ge 1$, it follows $f^{i-1}+f^{i}+f^{i+1}\ge 7$.

(b) $f^0+f^1\ge 5$ and $f^{n-1}+f^{n-2}\ge 5$.

Since $f^i\ge 2$ for every $i$ $(0\le i\le n-1)$, if $f^0\ge 3$, then $f^0+f^1\ge 5$; if $f^0=2$, without loss of generality, let $f(v_{0,0})=0$, $f(v_{0,1})=f(v_{0,2})=1$, then $f(v_{1,0})=2$, $f(v_{1,1})+f(v_{1,2})\ge 1$, i.e. $f^{1}\ge 3$, it follows $f^0+f^1\ge 5$.
Similarly, $f^{n-1}+f^{n-2}\ge 5$.

By (a) and (b), we have
\begin{equation*}\small{
\begin{split}
3w(f)&=\sum_{i=0}^{n-3}(f^{i}+f^{i+1}+f^{i+2})+2f^0+2f^{n-1}+f^1+f^{n-2}\\
     &=\sum_{i=0}^{n-3}(f^{i}+f^{i+1}+f^{i+2})+(f^0+f^1)+f^0+(f^{n-2}+f^{n-1})+f^{n-1}\\
     &\ge 7(n-2)+5+2+5+2=7n.
\end{split}}
\end{equation*}

Therefore, $\gamma_{oiR}(P_{n} \Box C_{3})\ge \lceil\frac{7n}{3}\rceil$.
\end{proof}

\section{Upper bounds on the outer independent Roman domination number of $P_{n}\Box C_{m}$, $n, m \geq 4$}

\begin{theorem}\label{thm:pncm}
For any integers $m, n\ge 4$, $\gamma_{oiR}(P_{n} \Box C_{m})\le \frac{5mn+5n+2m}{8}$.
\end{theorem}
\begin{proof}
Let $G=P_{n}\Box C_{m}$, $V(G)=\{v_{i,j}|0\leq i\leq n-1, 0\leq j\leq m-1\}$.
We first define a function $g$ on $P_{4} \Box C_{4}$ as follows.
\begin{eqnarray*}\small{
g(v_{i,j})=\left\{
\begin{array}{llllll}
2,   &(i=0, j=2)\vee (i=2, j=0),\\
1,   &(i=0, j=0)\vee (i=1,3, j=1,3)\vee (i=2, j=2),\\
0,   &otherwise.
\end{array}
\right. \text{\ \ i.e.\ \ }
g=\left(
\begin{array}{llllll}
1020 \\
0101 \\
2010 \\
0101
\end{array}
\right).}
\end{eqnarray*}

(1) For $m\equiv0(\bmod 4)$, we define an OIRDF $f$ on $P_{n} \Box C_{m}$ as the following.
\begin{equation}\small{
f(v_{i,j})=\left\{
\begin{array}{llllll}
2,                      &n\equiv0,2(\bmod\ 4), i=n-1, j\equiv3(\bmod\ 4),\\
g(v_{i\bmod4,j\bmod4}), &otherwise.
\end{array}\notag
\right.}
\end{equation}

Also, we can write $f$ as the following.
\begin{align*}\small{
\begin{array}{llllllllllll}
n\equiv1(\bmod\ 4),&&&n\equiv2(\bmod\ 4),\\
f=\left(\begin{array}{ccccccccccc}
1020&\cdots&1020\\
0101&\cdots&0101\\
2010&\cdots&2010\\
0101&\cdots&0101\\
\vdots&\ddots&\vdots\\
1020&\cdots&1020\\
0101&\cdots&0101\\
2010&\cdots&2010\\
0101&\cdots&0101\\
\hdashline
1020&\cdots&1020\\
\end{array}
\right),&&&
f=\left(
\begin{array}{cccccc}
1020&\cdots&1020\\
0101&\cdots&0101\\
2010&\cdots&2010\\
0101&\cdots&0101\\
\vdots&\ddots&\vdots\\
1020&\cdots&1020\\
0101&\cdots&0101\\
2010&\cdots&2010\\
0101&\cdots&0101\\
\hdashline
1020&\cdots&1020\\
0102&\cdots&0102
\end{array}\right),
\end{array}
}
\end{align*}
\begin{align*}\small{
\begin{array}{llllllllll}
n\equiv3(\bmod\ 4),&&&n\equiv0(\bmod\ 4),\\
f=\left(
\begin{array}{ccccccc}
1020&\cdots&1020\\
0101&\cdots&0101\\
2010&\cdots&2010\\
0101&\cdots&0101\\
\vdots&\ddots&\vdots\\
1020&\cdots&1020\\
0101&\cdots&0101\\
2010&\cdots&2010\\
0101&\cdots&0101\\
\hdashline
1020&\cdots&1020\\
0101&\cdots&0101\\
2010&\cdots&2010
\end{array}
\right),&&&
f=\left(
\begin{array}{cccccc}
1020&\cdots&1020\\
0101&\cdots&0101\\
2010&\cdots&2010\\
0101&\cdots&0101\\
\vdots&\ddots&\vdots\\
1020&\cdots&1020\\
0101&\cdots&0101\\
2010&\cdots&2010\\
0101&\cdots&0101\\
\hdashline
1020&\cdots&1020\\
0101&\cdots&0101\\
2010&\cdots&2010\\
0102&\cdots&0102\\
\end{array}\right).
\end{array}
}
\end{align*}

Then the weight is
\begin{equation}\small{
w(f)=\left\{
\begin{array}{lllllllll}
\frac{m}{4}\times \frac{n-4}{4}\times10+\frac{m}{4}\times11=\frac{5mn+2m}{8},     & n\equiv0(\bmod\ 4),\\
\frac{m}{4}\times \frac{n-1}{4}\times10+\frac{m}{4}\times3=\frac{5mn+m}{8},     & n\equiv1(\bmod\ 4),\\
\frac{m}{4}\times \frac{n-2}{4}\times10+\frac{m}{4}\times6=\frac{5mn+2m}{8},     & n\equiv2(\bmod\ 4),\\
\frac{m}{4}\times \frac{n-3}{4}\times10+\frac{m}{4}\times8=\frac{5mn+m}{8},     & n\equiv3(\bmod\ 4).
\end{array}\notag
\right.}
\end{equation}

Hence,
\begin{equation}
\gamma_{oiR}(P_n\Box C_m)\leq\left\{
\begin{array}{lllllllll}
\frac{5mn+2m}{8},      &n\equiv0,2(\bmod 4),\\
\frac{5mn+m}{8},       &n\equiv1,3(\bmod 4).
\end{array}\notag
\right.
\end{equation}

(2) For $m\equiv1(\bmod 4)$, we define an OIRDF $f$ on $P_{n} \Box C_{m}$ as the following.
\begin{equation}\small{
f(v_{i,j})=\left\{
\begin{array}{llllll}
2,                      &n\equiv0(\bmod\ 4), i=n-1, j\equiv3(\bmod\ 4),\\
                        &n\equiv2(\bmod\ 4), i=n-1, j\equiv3(\bmod\ 4), 0\le j\le m-5,\\
                        &n\equiv2(\bmod\ 4), i=n-1, j=m-1,\\
1,                      &i\equiv2(\bmod\ 4), j=m-2,\\
                        &n\equiv0,1,3(\bmod\ 4), i\equiv1,3(\bmod\ 4), j=m-1,\\
                        &n\equiv2(\bmod\ 4), i\equiv1(\bmod\ 4), j=m-1, 0\le i\le n-3,\\
                        &n\equiv2(\bmod\ 4), i\equiv3(\bmod\ 4), j=m-1,\\
0,                      &i\equiv2(\bmod\ 4), j=m-1,\\
g(v_{i\bmod4,j\bmod4}), &otherwise.
\end{array}\notag
\right.}
\end{equation}

Also, we can write $f$ as follows.
\begin{align*}\small{
\begin{array}{llllllllll}
n\equiv1(\bmod\ 4),&&&n\equiv2(\bmod\ 4),\\
f=\left(
\begin{array}{cccccccc}
1020&\cdots&1020&10201\\
0101&\cdots&0101&01011\\
2010&\cdots&2010&20110\\
0101&\cdots&0101&01011\\
\vdots&\ddots&\vdots&\vdots\\
1020&\cdots&1020&10201\\
0101&\cdots&0101&01011\\
2010&\cdots&2010&20110\\
0101&\cdots&0101&01011\\
\hdashline
1020&\cdots&1020&10201\\
\end{array}
\right),&&&
f=\left(
\begin{array}{cccccc}
1020&\cdots&1020&10201\\
0101&\cdots&0101&01011\\
2010&\cdots&2010&20110\\
0101&\cdots&0101&01011\\
\vdots&\ddots&\vdots&\vdots\\
1020&\cdots&1020&10201\\
0101&\cdots&0101&01011\\
2010&\cdots&2010&20110\\
0101&\cdots&0101&01011\\
\hdashline
1020&\cdots&1020&10201\\
0102&\cdots&0102&01012\\
\end{array}\right),
\end{array}
}
\end{align*}
\begin{align*}\small{
\begin{array}{lllllllllll}
n\equiv3(\bmod\ 4),&&&n\equiv0(\bmod\ 4),\\
f=\left(
\begin{array}{cccccccc}
1020&\cdots&1020&10201\\
0101&\cdots&0101&01011\\
2010&\cdots&2010&20110\\
0101&\cdots&0101&01011\\
\vdots&\ddots&\vdots&\vdots\\
1020&\cdots&1020&10201\\
0101&\cdots&0101&01011\\
2010&\cdots&2010&20110\\
0101&\cdots&0101&01011\\
\hdashline
1020&\cdots&1020&10201\\
0101&\cdots&0101&01011\\
2010&\cdots&2010&20110\\
\end{array}
\right),&&&
\  f=\left(
\begin{array}{cccccc}
1020&\cdots&1020&10201\\
0101&\cdots&0101&01011\\
2010&\cdots&2010&20110\\
0101&\cdots&0101&01011\\
\vdots&\ddots&\vdots&\vdots\\
1020&\cdots&1020&10201\\
0101&\cdots&0101&01011\\
2010&\cdots&2010&20110\\
0101&\cdots&0101&01011\\
\hdashline
1020&\cdots&1020&10201\\
0101&\cdots&0101&01011\\
2010&\cdots&2010&20110\\
0102&\cdots&0102&01021\\
\end{array}\right).
\end{array}
}
\end{align*}

Then the weight is
\begin{equation}\small{
w(f)=\left\{
\begin{array}{lllllllll}
\frac{m-5}{4}\times \frac{n-4}{4}\times10+\frac{n-4}{4}\times14+\frac{m-5}{4}\times11+15=\frac{5mn+3n+2m-2}{8},     & n\equiv0(\bmod\ 4),\\
\frac{m-5}{4}\times \frac{n-1}{4}\times10+\frac{n-1}{4}\times14+\frac{m-5}{4}\times3+4=\frac{5mn+3n+m-1}{8},     & n\equiv1(\bmod\ 4),\\
\frac{m-5}{4}\times \frac{n-2}{4}\times10+\frac{n-2}{4}\times14+\frac{m-5}{4}\times6+8=\frac{5mn+3n+2m-2}{8},     & n\equiv2(\bmod\ 4),\\
\frac{m-5}{4}\times \frac{n-3}{4}\times10+\frac{n-3}{4}\times14+\frac{m-5}{4}\times8+11=\frac{5mn+3n+m-1}{8},    & n\equiv3(\bmod\ 4).\\
\end{array}\notag
\right.}
\end{equation}

Hence,
\begin{equation}
\gamma_{oiR}(P_n\Box C_m)\leq\left\{
\begin{array}{lllllllll}
\frac{5mn+3n+2m-2}{8},      &n\equiv0,2(\bmod 4),\\
\frac{5mn+3n+m-1}{8},       &n\equiv1,3(\bmod 4).
\end{array}\notag
\right.
\end{equation}

(3) For $m\equiv2(\bmod 4)$, we define an OIRDF $f$ on $P_{n} \Box C_{m}$ as the following.
\begin{equation}\small{
f(v_{i,j})=\left\{
\begin{array}{llllll}
2,                      &i\equiv3(\bmod\ 4), j=m-3,\\
                        &i\equiv1(\bmod\ 4), j=m-1,\\
                        &n\equiv0(\bmod\ 4), i=n-1, j\equiv3(\bmod\ 4),\\
                        &n\equiv2(\bmod\ 4), i=n-1, j\equiv3(\bmod\ 4), 0\le j\le m-6,\\
1,                      &n\equiv0,1,2(\bmod\ 4),i\equiv2(\bmod\ 4), j=m-2,\\
                        &n\equiv3(\bmod\ 4),i\equiv2(\bmod\ 4), j=m-2, 0\le i\le n-4,\\
                        &n\equiv1(\bmod\ 4), i=n-1, j=m-1,\\
g(v_{i\bmod4,j\bmod4}), &otherwise.
\end{array}\notag
\right.}
\end{equation}

Also, we can write $f$ as follows.
\begin{align*}\small{
\begin{array}{llllllllll}
n\equiv1(\bmod\ 4),&&&n\equiv2(\bmod\ 4),\\
 f=\left(
\begin{array}{cccccccc}
1020&\cdots&1020&102010\\
0101&\cdots&0101&010102\\
2010&\cdots&2010&201010\\
0101&\cdots&0101&010201\\
\vdots&\ddots&\vdots&\vdots\\
1020&\cdots&1020&102010\\
0101&\cdots&0101&010102\\
2010&\cdots&2010&201010\\
0101&\cdots&0101&010201\\
\hdashline
1020&\cdots&1020&102011\\
\end{array}
\right),&&&
f=\left(
\begin{array}{cccccc}
1020&\cdots&1020&102010\\
0101&\cdots&0101&010102\\
2010&\cdots&2010&201010\\
0101&\cdots&0101&010201\\
\vdots&\ddots&\vdots&\vdots\\
1020&\cdots&1020&102010\\
0101&\cdots&0101&010102\\
2010&\cdots&2010&201010\\
0101&\cdots&0101&010201\\
\hdashline
1020&\cdots&1020&102010\\
0102&\cdots&0102&010102\\
\end{array}
\right),
\end{array}
}
\end{align*}
\begin{align*}\small{
\begin{array}{llllllllllll}
n\equiv3(\bmod\ 4),&&& n\equiv0(\bmod\ 4),\\
f=\left(
\begin{array}{cccccccc}
1020&\cdots&1020&102010\\
0101&\cdots&0101&010102\\
2010&\cdots&2010&201010\\
0101&\cdots&0101&010201\\
\vdots&\ddots&\vdots&\vdots\\
1020&\cdots&1020&102010\\
0101&\cdots&0101&010102\\
2010&\cdots&2010&201010\\
0101&\cdots&0101&010201\\
\hdashline
1020&\cdots&1020&102010\\
0101&\cdots&0101&010102\\
2010&\cdots&2010&201020\\
\end{array}
\right),&&&
f=\left(
\begin{array}{cccccc}
1020&\cdots&1020&102010\\
0101&\cdots&0101&010102\\
2010&\cdots&2010&201010\\
0101&\cdots&0101&010201\\
\vdots&\ddots&\vdots&\vdots\\
1020&\cdots&1020&102010\\
0101&\cdots&0101&010102\\
2010&\cdots&2010&201010\\
0101&\cdots&0101&010201\\
\hdashline
1020&\cdots&1020&102010\\
0101&\cdots&0101&010102\\
2010&\cdots&2010&201010\\
0102&\cdots&0102&010201\\
\end{array}
\right).
\end{array}
}
\end{align*}

Then the weight is
\begin{equation}\small{
w(f)=\left\{
\begin{array}{lllllllll}
\frac{m-6}{4}\times \frac{n-4}{4}\times10+\frac{n-4}{4}\times16+\frac{m-6}{4}\times11+16=\frac{5mn+2n+2m-12}{8},     & n\equiv0(\bmod\ 4),\\
\frac{m-6}{4}\times \frac{n-1}{4}\times10+\frac{n-1}{4}\times16+\frac{m-6}{4}\times3+5=\frac{5mn+2n+m+2}{8},     & n\equiv1(\bmod\ 4),\\
\frac{m-6}{4}\times \frac{n-2}{4}\times10+\frac{n-2}{4}\times16+\frac{m-6}{4}\times6+8=\frac{5mn+2n+2m-12}{8},     & n\equiv2(\bmod\ 4),\\
\frac{m-6}{4}\times \frac{n-3}{4}\times10+\frac{n-3}{4}\times16+\frac{m-6}{4}\times8+13=\frac{5mn+2n+m+2}{8},    & n\equiv3(\bmod\ 4).\\
\end{array}\notag
\right.}
\end{equation}

Hence,
\begin{equation}\small{
\gamma_{oiR}(P_n\Box C_m)\leq\left\{
\begin{array}{lllllllll}
\frac{5mn+2n+2m-12}{8},     &n\equiv0,2(\bmod 4),\\
\frac{5mn+2n+m+2}{8},       &n\equiv1,3(\bmod 4).
\end{array}\notag
\right.}
\end{equation}

(4) For $m\equiv3(\bmod 4)$, we define an OIRDF $f$ on $P_{n} \Box C_{m}$ as the following.
\begin{equation}\small{
f(v_{i,j})=\left\{
\begin{array}{llllll}
2,                      &i\equiv1(\bmod\ 4), j=m-2,\\
                        &n\equiv0(\bmod\ 4), i=n-1, j\equiv3(\bmod\ 4),\\
1,                      &i\equiv3(\bmod\ 4), j=m-1,\\
                        &n\equiv0,3(\bmod\ 4),i\equiv0,1(\bmod\ 4), j=m-1,\\
                        &n\equiv0(\bmod\ 4),i=n-2, j=m-3,\\
                        &n\equiv1(\bmod\ 4),i\equiv0(\bmod\ 4), j=m-1, 0\le i\le n-5,\\
                        &n\equiv1(\bmod\ 4),i\equiv1(\bmod\ 4), j=m-1,\\
                        &n\equiv2(\bmod\ 4),i\equiv0(\bmod\ 4), j=m-1,\\
                        &n\equiv2(\bmod\ 4),i\equiv1(\bmod\ 4), j=m-1, 0\le i\le n-5,\\
                        &n\equiv2(\bmod\ 4), i=n-1, j\equiv0(\bmod\ 4), 0\le j\le m-7,\\
g(v_{i\bmod4,j\bmod4}), &otherwise.
\end{array}\notag
\right.}
\end{equation}

Also, we can write $f$ as follows.
\begin{align*}\small{
\begin{array}{lllllllll}
n\equiv1(\bmod\ 4),&&&n\equiv2(\bmod\ 4),\\
f=\left(
\begin{array}{cccccccc}
1020&\cdots&1020&101\\
0101&\cdots&0101&021\\
2010&\cdots&2010&201\\
0101&\cdots&0101&011\\
\vdots&\ddots&\vdots&\vdots\\
1020&\cdots&1020&101\\
0101&\cdots&0101&021\\
2010&\cdots&2010&201\\
0101&\cdots&0101&011\\
\hdashline
1020&\cdots&1020&102\\
\end{array}
\right),&&&
f=\left(
\begin{array}{cccccc}
1020&\cdots&1020&101\\
0101&\cdots&0101&021\\
2010&\cdots&2010&201\\
0101&\cdots&0101&011\\
\vdots&\ddots&\vdots&\vdots\\
1020&\cdots&1020&101\\
0101&\cdots&0101&021\\
2010&\cdots&2010&201\\
0101&\cdots&0101&011\\
\hdashline
1020&\cdots&1020&101\\
1101&\cdots&1101&020\\
\end{array}
\right),
\end{array}
}
\end{align*}
\begin{align*}\small{
\begin{array}{llllllllll}
n\equiv3(\bmod\ 4),& & &n\equiv0(\bmod\ 4),\\
f=\left(
\begin{array}{cccccccc}
1020&\cdots&1020&101\\
0101&\cdots&0101&021\\
2010&\cdots&2010&201\\
0101&\cdots&0101&011\\
\vdots&\ddots&\vdots&\vdots\\
1020&\cdots&1020&101\\
0101&\cdots&0101&021\\
2010&\cdots&2010&201\\
0101&\cdots&0101&011\\
\hdashline
1020&\cdots&1020&101\\
0101&\cdots&0101&021\\
2010&\cdots&2010&201\\
\end{array}
\right),& & &
f=\left(
\begin{array}{cccccc}
1020&\cdots&1020&101\\
0101&\cdots&0101&021\\
2010&\cdots&2010&201\\
0101&\cdots&0101&011\\
\vdots&\ddots&\vdots&\vdots\\
1020&\cdots&1020&101\\
0101&\cdots&0101&021\\
2010&\cdots&2010&201\\
0101&\cdots&0101&011\\
\hdashline
1020&\cdots&1020&101\\
0101&\cdots&0101&021\\
2010&\cdots&2010&101\\
0102&\cdots&0102&011\\
\end{array}
\right).
\end{array}
}
\end{align*}

Then the weight is
\begin{equation}\small{
w(f)=\left\{
\begin{array}{lllllllll}
\frac{m-3}{4}\times \frac{n-4}{4}\times10+\frac{n-4}{4}\times10+\frac{m-3}{4}\times11+9=\frac{5mn+5n+2m-14}{8},     & n\equiv0(\bmod\ 4),\\
\frac{m-3}{4}\times \frac{n-1}{4}\times10+\frac{n-1}{4}\times10+\frac{m-3}{4}\times3+3=\frac{5mn+5n+m+1}{8},       & n\equiv1(\bmod\ 4),\\
\frac{m-3}{4}\times \frac{n-2}{4}\times10+\frac{n-2}{4}\times10+\frac{m-3}{4}\times6+4=\frac{5mn+5n+2m-14}{8},     & n\equiv2(\bmod\ 4),\\
\frac{m-3}{4}\times \frac{n-3}{4}\times10+\frac{n-3}{4}\times10+\frac{m-3}{4}\times8+8=\frac{5mn+5n+m+1}{8},        & n\equiv3(\bmod\ 4).\\
\end{array}\notag
\right.}
\end{equation}

Hence,
\begin{equation}\small{
\gamma_{oiR}(P_n\Box C_m)\leq\left\{
\begin{array}{lllllllll}
\frac{5mn+5n+2m-14}{8},       &n\equiv0,2(\bmod 4),\\
\frac{5mn+5n+m+1}{8},       &n\equiv1,3(\bmod 4).\\
\end{array}\notag
\right.}
\end{equation}

Therefore, $\gamma_{oiR}(P_{n} \Box C_{m})\le \frac{5mn+5n+2m}{8}$.
\end{proof}

\begin{acknowledgements}
The authors would like to thank the anonymous referee, whose suggestions greatly
improved the exposition of this paper.
\end{acknowledgements}

%

\nocite{*}
\bibliographystyle{fundam}
\bibliography{citations}


\end{document}